\documentclass[11pt]{amsart}
\usepackage{amssymb,amsmath,amsthm,amscd,mathrsfs,graphicx}
\usepackage[cmtip,all]{xy}
\usepackage{pb-diagram}
\usepackage{tikz}

\numberwithin{equation}{section}
\newtheorem{teo}{Theorem}[section]

\newtheorem{theorem}{Theorem}[section]

\newtheorem{question}{Question}[section]

\newcommand{\ov}[1]{\overline{#1}}

\newcommand{\ve}{\varepsilon}

\theoremstyle{definition}

\newtheorem{dfn}[teo]{Definition}

\theoremstyle{remark}
\newtheorem{remark}[teo]{Remark}

\begin{document}
\bibliographystyle{amsplain}

\title[Concavity of solutions to semilinear equations]{Concavity of solutions to semilinear \\ equations in dimension two}

\author[A. Chau]{Albert Chau}
\address{Department of Mathematics, The University of British Columbia, 1984 Mathematics Road, Vancouver, B.C.,  Canada V6T 1Z2.  Email: chau@math.ubc.ca. } 
\author[B. Weinkove]{Ben Weinkove}
\address{Department of Mathematics, Northwestern University, 2033 Sheridan Road, Evanston, IL 60208, USA.  Email: weinkove@math.northwestern.edu.}

\thanks{Research supported in part by  NSERC grant $\#$327637-06 and NSF grant DMS-2005311.}

\maketitle

\vspace{-20pt}

\begin{abstract}  We consider the Dirichlet problem for a class of semilinear equations on two dimensional convex domains.  We give a sufficient condition for the solution to be concave.  Our condition uses comparison with ellipses, and is motivated by an idea of Kosmodem'yanskii.

We also prove a result on propagation of concavity of solutions from the boundary, which holds in all dimensions.
\end{abstract}

\section{Introduction}

Let $\Omega \subset \mathbb{R}^2$ be a bounded convex domain with smooth boundary and let $u \in C^{\infty}(\ov{\Omega})$ solve the semilinear equation
\begin{equation}\label{mainequ}
\begin{split}
-\Delta u = {} & f(u), \qquad \textrm{in } \Omega, \\
u = {} & 0, \qquad \textrm{on } \partial \Omega.
\end{split}
\end{equation}
for a smooth positive function $f: \mathbb{R} \rightarrow \mathbb{R}$.

In the special case $f=1$,  the corresponding Dirichlet problem,
\begin{equation}\label{torsion}
\begin{split}
-\Delta u = {} & 1, \qquad \textrm{in } \Omega, \\
u = {} & 0, \qquad \textrm{on } \partial \Omega.
\end{split}
\end{equation}
is known as the \emph{torsion problem}.

A natural question is whether the solution $u$ should inherit concavity properties from the domain $\Omega$.   In general $u$ is not concave, even in the case of the torsion problem.  Indeed  the solution $u$ to (\ref{torsion}) on an equilateral triangle is not concave \cite{KM}.  One can obtain a non-concave solution with smooth boundary by taking $\Omega = \{ u\ge \ve\}$ for a small $\ve>0$.

Recently Steinerberger \cite{S} asked:

\begin{question} Under what conditions on $f$ and $\Omega$ is the function $u$ concave on $\Omega$? \end{question}

There appear to be very few results on this.  For the torsion problem (\ref{torsion}), Kosmodem'yanskii \cite{K} showed that $u$ is concave if $\Omega$ satisfies a third order contact parabola condition, which we discuss below in Section \ref{comparison}.

  For the torsion problem, the answer is clearly yes when $\Omega$  is the interior of an ellipse, say, given by $a^2x^2 + b^2 y^2=1$.  Indeed  the solution  is the concave quadratic function
$$u = c(1- a^2x^2 - b^2 y^2), \quad \textrm{for } c =\left( 2a^2 + 2b^2 \right)^{-1}.$$
Solutions will also be concave if $\Omega$ is a small perturbation of an ellipse (see \cite{HNST}, for example).  On the other hand, as discussed above, $u$ is not concave for some domains where $\partial \Omega$ is close to a triangle.
It is then natural to ask whether concavity for $u$ holds when $\Omega$ is ``not too far'' from being an ellipse.   

Our main result gives a quantitative sufficient condition on $\Omega$ for concavity 
 for a large class of $f$ which includes the case $f=1$.  Our condition uses comparison with ellipses.

\pagebreak[3]

\begin{theorem} \label{mainthm}  Let $u$ be a solution of (\ref{mainequ}).  Assume:
 \begin{enumerate}
 	\item For each point $p \in \partial \Omega$, there is an ellipse $E$, enclosing a finite area $A$, containing $\Omega$ and tangent to $\partial \Omega$ at $p$ such that, at $p$, 
	\begin{equation} \label{assume1}
	K_E \left( \frac{\pi^{2/3}}{ A^{2/3}K_E^{4/3}} + 1 +  \frac{(\partial_s K_{E})^2}{9K_E^4} \right) \ge K \left( 1+ \frac{(\partial_sK)^2}{9K^4} \right),
	\end{equation}
	where $K$ and $K_E$ are the curvatures of $\partial \Omega$ and $E$ respectively, and $\partial_s$ is the derivative with respect to arclength.  
		\item $f(0)=1$, and on $[0, \max_{\ov{\Omega}} u]$ we have $f>0, \ f''\le 0$ and $f'\le 0$. 
\end{enumerate}
Then $u$ is concave on $\overline \Omega$.
 \end{theorem}

\begin{remark}\label{remark}
\ 
\begin{enumerate}

\item[(i)]  Condition (1) above is clearly satisfied when $\partial \Omega$ is itself an ellipse.  Indeed, given any $p\in \partial \Omega$ we simply take $E=\partial \Omega$ in (1) in which case we replace $K_E, \partial_s K_E$ with $K, \partial_s K$ in \eqref{assume1}, leaving a strict inequality that is obviously satisfied since $\pi^{2/3}/(A^{2/3} K^{4/3})>0$. In particular, condition (1) will
continue to hold for all domains whose $C^3$ distance to a given ellipse is
comparable to the term $\pi^{2/3}/(A^{2/3}K^{4/3})$.

\smallskip

\item[(ii)] Our proof is motivated by the argument of Kosmodem'yanskii \cite{K} who used a comparison with parabolas which have a third order contact with the boundary of $\Omega$.   
 In Section \ref{comparison} below, we describe a key difference between these two conditions:  unlike (1) above, Kosmodem'yanskii's condition cannot hold in any $C^3$ neighborhood of any domain.  In this way one can produce domains that satisfy our conditions but not those of \cite{K}, see Remark \ref{remark2} below.  
 
 \smallskip
 
 \item[(iii)] Examples of functions $f$ satisfying (2) and for which there exist smooth solutions to (\ref{mainequ}) are $f(u) = 1-cu^k$ for $k=1,2,3, \ldots$ and $c$ a sufficiently small positive constant, depending only on $\Omega$ and $k$.

 \smallskip
 
\item[(iv)] When $f=1$ the conclusion can be strengthened to $u$ being strongly concave, namely $D^2u <0$ on $\overline{\Omega}$ (see Remark \ref{remarksection2} (ii) below).  In particular, this implies that the conditions in Theorem \ref{mainthm} are not optimal in general.

\smallskip

\item[(v)]  It's not clear whether the analogue of Theorem \ref{mainthm} holds for dimensions larger than 2.  Using our approach,   new difficulties arise from the third order derivatives of $u$.

 \end{enumerate}
  
\end{remark}

By our methods we also obtain a short proof of the following   ``propagating concavity from the boundary'' result, which holds in all dimensions. 

\begin{theorem}\label{thmgeneralstbg} Let $\Omega \subset \mathbb{R}^n$ be a bounded  convex domain with smooth boundary and let $u \in C^{\infty}(\ov{\Omega})$ solve (\ref{mainequ}).  
Assume $f>0$ and $f''\le 0$.  
If $D^2u \le 0$ on $\partial \Omega$ then $u$ is concave in $\Omega$.
\end{theorem}

This extends a recent result of Steinerberger \cite{S}, who proved the same result under additional assumptions on $f'$.  Earlier work of Keady-McNabb \cite{KM} established it for the torsion problem.

We end the introduction with some historical remarks.  Instead of asking whether $u$ is concave, one could ask  whether weaker concavity properties hold.  Lions \cite{L} conjectured in 1981 that for a general $f$ the superlevel sets of $u$ should be convex (i.e. $u$ is \emph{quasiconcave}), and this has led to many further investigations.  By a classical result of Gidas-Ni-Nirenberg \cite{GNN}, $u$ is quasiconcave  in the special case when $\Omega$ is a disc, since then $u$ is radially symmetric.   Makar-Limanov showed that solutions of the two-dimension torsion problem (\ref{torsion}) are quasiconcave since $\sqrt{u}$ is concave \cite{M}.  In a well-known paper, Brascamp and Lieb \cite{BL} showed that if $f(u)= \lambda u$ for $\lambda>0$ then $\log u$ is concave.  There have been many further extensions and generalizations.  We refer the reader to, for example,  \cite{APP, ALL, BGMX, CF, CS, GP, Kaw, Kea, Ken, Kor, KL} and the references therein.   On the other hand, a 2016 counterexample of Hamel-Nadirashvili-Sire \cite{HNS} shows that quasiconcavity is false for general $f$.

The outline of this paper is as follows.  In Sections \ref{section2} and \ref{section3} we prove Theorems \ref{mainthm} and \ref{thmgeneralstbg} respectively.  In Section \ref{comparison} we compare our condition (1) in Theorem \ref{mainthm} to the parabola condition of Kosmodem'yanskii \cite{K}. 
 	
\section{Proof of Theorem \ref{mainthm}} \label{section2}
 	
Assume for a contradiction that $u$ is not concave on $\overline{\Omega}$.  Then there exists a  
 vector $e=(e_x, e_y)$ and a point $p \in \ov{\Omega}$ such that the quantity \begin{equation*}\label{ueeeeee} H=u_{ee}  \end{equation*} 
 achieves a strictly positive maximum at $p$.  Here and henceforth we are using subscripts to denote partial derivatives.

We first show that we may assume that $p$ lies in the boundary $\partial \Omega$.  Consider the set
$$M =\{ x \in \Omega \ | \ H(x) = \sup_{\ov{\Omega}} H \}.$$
This is a closed set in $\Omega$.  But it is also open, since if $x \in M$ then in a small ball $B$ with $x\in B \subset \Omega$ we have $$\Delta H= \Delta u_{ee} =-f'(u)u_{ee} -f''(u) u_e^2 \ge 0,$$
since $f' \le 0$ and $f''\leq 0$.  The strong maximum principle implies that $H$ is constant on $B$ and hence $B \subset \Omega$.  Hence $M$ is closed and open in $\Omega$, so $M=\Omega$ or $M$ is empty.  Either way, $H$ achieves its maximum at a boundary point $p \in \partial \Omega$.

  We may rotate and translate the coordinates so that $p$ is the origin and $\partial \Omega$ is locally given by $y=\rho(x)$ for $\rho$ a 	convex function with $\rho(0)=0$ and $\rho'(0)=0$.  Assumption (1) of the theorem implies that $\rho''(0)>0$.
	
	\begin{dfn} \label{preferred}
	We say that such a coordinate system  is a \emph{preferred coordinate system} at $p\in \partial \Omega$.  
 	 \end{dfn}
	 
 We make the following claim.
 
\bigskip
\noindent
{\bf Claim.} 
At $p$, we have
 \begin{equation}\label{uxxininininin}
 	0<u_y \le \frac{9(\rho'')^3}{9(\rho'')^4 + (\rho''')^2}.
 \end{equation} 
where the second inequality is an equality if and only if $u$ is quadratic and $\partial \Omega$ is an ellipse.

\begin{proof}[Proof of Claim] The first inequality is an immediate consequence of the maximum principle and the Hopf Lemma.  Indeed since $\Delta u =-f(u) < 0$, we have $u>0$ in $\Omega$ and $\partial_{\nu} u<0$ on $\partial \Omega$, where $\nu$ is the \emph{outward} pointing normal.

We will prove the second inequality using conditions (1) and (2) of Theorem \ref{mainthm}.  We recall a basic fact that any ellipse in the plane is given by an equation
 \begin{equation} \label{generalellipse}
 a^2 (x-h)^2 + b(x-h)(y-k)+c^2(y-k)^2=1,
 \end{equation}
 for constants $a, b, c, h$ and $k$ with $D:= 4a^2 c^2 - b^2>0$.   Moreover, if $A$ is the area of the ellipse then 
 \begin{equation} \label{DA}
 D = \frac{4\pi^2}{A^2}.
 \end{equation}
 Let $E$ be the ellipse as in condition (1) of the theorem.
  We may and do assume without loss of generality that $a$ and $c$ are positive.  
  
  By the assumptions on $E$, the origin is a point on the ellipse and the tangent line to $E$ at the origin coincides with the $x$-axis.  Hence if we express $E$ near the origin as a graph $y=g(x)$ then $g(0)=0$ and $g'(0)=0$.  Since the ellipse contains $\Omega$, which is in the upper half plane, we have $g''(0)>0$.  
 
Evaluating (\ref{generalellipse}), and its implicit derivative with respect to $x$, at the origin, we can solve for $h$, $k$ to obtain
$$h=-\frac{b}{a\sqrt{D}}, \quad k = \frac{2a}{\sqrt{D}}.$$
Differentiating (\ref{generalellipse}) implicitly another two times and evaluating at the origin, we obtain
\begin{equation} \label{g}
g''(0) = \frac{2a^3}{\sqrt{D}}, \ g'''(0) = \frac{6a^4b}{D}.
\end{equation}

Define function $w=w(x,y)$ by
\[
\begin{split}
	w= {} & \frac{1}{2(a^2+c^2)} \bigg( 1 - a^2 \left( x+ \frac{b}{a\sqrt{D}} \right)^2 \\ {} &  - b\left( x + \frac{b}{a\sqrt{D}} \right) \left( y- \frac{2a}{\sqrt{D}} \right)  - c^2 \left( y -\frac{2a}{\sqrt{D}}\right)^2 \bigg).
\end{split}
\] 
Then $\Delta w=-1$ and the set $\{ w>0 \}$ coincides with the interior of the ellipse $E$.  By assumption the ellipse contains the set $\Omega$.  Hence 
$$\Delta (w-u)=f(u)-1 \le 0, \quad \textrm{on } \Omega,$$ and $w-u \ge 0$ on $\partial \Omega$, where  we used the condition that $f' \le 0$.
Moreover $w-u=0$ at the origin.  The strong maximum principle then implies that either $w-u =0$ on $\Omega$ in which case $u=w$ is a quadratic and $\partial \Omega$ is the ellipse $E$, or else $w-u>0$ on  $\Omega$ and $w-u=0$ at $p\in \partial \Omega$.  Hence 
\begin{equation} \label{uywy}
u_y(0,0) \le w_y(0,0) = \frac{\sqrt{D}}{2a(a^2+c^2)}.
\end{equation}
By the Hopf Lemma, the first inequality of (\ref{uywy}) is strict unless $u$ is quadratic.
We now wish to write the right hand sides of (\ref{uxxininininin}) and (\ref{uywy}) in terms of the intrinsic quantities $K, K_E, \partial_s K, \partial_s K_E$ and $A$  featured in (\ref{assume1}).
Note that, evaluating at $p$, using (\ref{g}),
\begin{equation} \label{KE}
K_E = g''(0) = \frac{2a^3}{\sqrt{D}}, \ \partial_sK_E = g'''(0) = \frac{6a^4b}{D}.
\end{equation}
Recalling the formula $D=4a^2c^2 - b^2$, we have
\begin{equation}  \label{uno}
\begin{split}
\frac{\sqrt{D}}{2a(a^2+c^2)} =  {} & \frac{\sqrt{D}}{2a^3 + (2a)^{-1} (D+b^2)} \\
= {} &  K_E^{-1} \left( \frac{\pi^{2/3}}{K_E^{4/3} A^{2/3}} + 1 + \frac{(\partial_s K_E)^2}{9K_E^4} \right)^{-1},
\end{split}
\end{equation}
after substituting for $D$, $a$ and $b$ from (\ref{DA}) and (\ref{KE}), and simplifying. 
On the other hand, at $p$ we have
$$K = \rho''(0), \ \partial_s K = \rho'''(0),$$
and hence
\begin{equation} \label{dos}
\frac{9(\rho'')^3}{9(\rho'')^4 + (\rho''')^2} = K^{-1} \left( 1+ \frac{(\partial_s K)^2}{9 K^4} \right)^{-1}.
\end{equation} 
Then (\ref{uxxininininin}) follows from combining (\ref{uywy}), (\ref{uno}), (\ref{dos}) and (\ref{assume1}), completing the proof of the claim. \end{proof}

 		Differentiating the equation $u(x,\rho(x))=0$  and evaluating at $x=0$, we obtain, recalling that $\rho'(0)=0$,
 	\begin{equation} \label{dom}
 		\begin{split}
 			u_x= {} & 0 \\
 			u_{xx} + u_y \rho'' = {} & 0 \\
 			u_{xxx} + 3u_{xy} \rho'' + u_y \rho'''= {} & 0.
 		\end{split}
 	\end{equation}

We recall that $e$ is the unit vector used in the definition of $H$. If the vector $e$ is proportional to the vector $(1,0)$ then we have $u_{xx} > 0$ at $p$.  But since $u_y>0$ and $\rho''(0)>0$ this contradicts the second equation of (\ref{dom}).  
 
If $e$ is proportional to the vector $(0,1)$ then $u_{yy} > 0$ at $p$.  But then at $p$ we have $$u_{xx} = -f(u)-u_{yy} < -1$$ from the equation $\Delta u =-f(u)$ and the fact that $f(0)=1$.  But by (\ref{uxxininininin}) and  the second equation of (\ref{dom}) we have
 $$u_{xx} = - u_y \rho'' \ge - \frac{9(\rho'')^4}{9(\rho'')^4+ (\rho''')^2} \ge -1,$$
 a contradiction.
 
 Hence we may assume, after rescaling, that $e=(1,\tau)$ for a real nonzero number $\tau$.  Write
 \begin{equation}\label{H}
 H = u_{ee} = u_{xx}  + \tau^2 u_{yy} + 2\tau u_{xy}.
 \end{equation}

 	Note that at $p$ we have $u_x=0$ and hence, differentiating the equation $u_{xx}+u_{yy}= - f(u)$ we get
 	\begin{equation} \label{diffeqn}
 		\begin{split} 
 			u_{xxx} + u_{yyx} = {} & 0, \\
 			u_{xxy} + u_{yyy} = {} & - f'(u) u_y.
 		\end{split}
 	\end{equation}
 	Then, recalling that $H$ achieves a maximum at $p$,
 	\[
 	\begin{split}
 		0 = H_x = {} & u_{eex} \\
 		= {} & u_{xxx} + \tau^2 u_{yyx} + 2\tau u_{xyx} \\
 		= {} & (1-\tau^2) u_{xxx} - 2\tau u_{yyy} - 2\tau f'(u) u_y.
 	\end{split}
 	\]
 	Hence
 	$$u_{yyy} = \frac{1-\tau^2}{2\tau} u_{xxx}- f'(u)u_y.$$
 	Next, using this and \eqref{diffeqn}, 
 	\[
 	\begin{split}
 		0\ge H_y = {} & u_{xxy} + \tau^2 u_{yyy} + 2\tau u_{xyy} \\
 		= {} & (\tau^2-1) u_{yyy} - 2\tau u_{xxx}- f'(u) u_y \\
 		= {} & - \frac{(1-\tau^2)^2 + 4\tau^2}{2\tau} u_{xxx} - (\tau^2-1) f'(u) u_y - f'(u)u_y\\
 		= {} & - \frac{(1+\tau^2)^2}{2\tau} u_{xxx} - \tau^2 f'(u) u_y.
 	\end{split}
 	\]
 	Hence 
 	$$\tau u_{xxx} \ge - \frac{2\tau^4}{(1+\tau^2)^2} f'(u)u_y.$$
 	Compute using  (\ref{H}) and (\ref{dom}),
 	\[
 	\begin{split}
 		-\frac{2\tau^4}{(1+\tau^2)^2} f'(u)u_y \le \tau u_{xxx} = {} & - 3\tau u_{xy}\rho'' - \tau u_y \rho''' \\
 		= {} & \frac{3\rho''}{2} (-H +  u_{xx} + \tau^2 u_{yy}) - \tau u_y \rho'''. 
 	\end{split}
 	\]
 	Using $u_{xx}+u_{yy}=-1$ at $p$ (using $f(0)=1$) and $u_{xx} = -u_y \rho''$ we obtain, recalling that $H$ is positive at $p$,
 	\begin{equation} \label{buenof}
 		\begin{split}
 			0 <  3\rho'' H 
 			\le {} & 3\rho'' (u_{xx} + \tau^2(-1-u_{xx})) - 2\tau u_y \rho''' + \frac{4\tau^2}{(1+\tau^2)^2} f'(u) u_y\\
 			= {} & 3\rho'' (-u_y \rho'' + \tau^2 (-1+ u_y \rho'')) - 2\tau u_y \rho'''+ \frac{4\tau^2}{(1+\tau^2)^2} f'(u) u_y\\
 			\le {} & -u_y \left\{ \left(  \frac{3\rho''}{u_y} - 3 (\rho'')^2 \right) \tau^2 + 2\rho'''\tau + 3(\rho'')^2  \right\},
 		\end{split}
 	\end{equation}
 	where for the last line we used $f' \le 0$ and $u_y > 0$.
 	
 	 	The inequality (\ref{buenof}) implies that the quadratic function
 	$$q(\tau) = \left( \frac{3\rho''}{u_y} - 3 (\rho'')^2 \right) \tau^2 + 2\rho'''\tau + 3(\rho'')^2,$$
 	is negative for some $\tau$ and hence has two roots.   Here we are using the fact that from (\ref{uxxininininin}) we have
 	$$ \frac{3\rho''}{u_y} - 3 (\rho'')^2 = \frac{3\rho''}{u_y} (  1 - \rho'' u_y) > 0.$$
	where we have assumed without loss of generality that $u$ is not quadratic, in which case the theorem would be trivial.

	 Since $q$ has two roots, we have 	$$4(\rho''')^2 - 4 \left( \frac{3\rho''}{u_y} - 3 (\rho'')^2 \right) 3(\rho'')^2 > 0,$$
 	namely
 	$$(\rho''')^2 + 9 (\rho'')^4 >  \frac{9(\rho'')^3}{u_y},$$
 	which contradicts the claim.  This completes the proof of the theorem.

\setcounter{teo}{0}

\begin{remark} \ \, \label{remarksection2}
\begin{enumerate}
\item[(i)]  Kosmodem'yanskii \cite{K} uses some similar arguments to those above.  He compares using a parabola instead of an ellipse, and applies the maximum principle to the determinant of the Hessian of $u$ rather than $u_{ee}$ for some fixed vector $e$.
\item[(ii)] In the case of the torsion problem (\ref{torsion}), one can strengthen the conclusion of Theorem \ref{mainthm} to $u$ being \emph{strongly convex}, namely that $D^2 u<0$ on $\overline{\Omega}$.  Indeed, at the beginning of the proof, one can assume for a contradiction that $H=u_{ee}$ has a nonnegative (rather than strictly positive) maximum.  Since $\Delta H=0$, this must occur at a boundary point.  The rest of the proof goes through in the same way.
\end{enumerate}
\end{remark}

\section{Proof of Theorem \ref{thmgeneralstbg}} \label{section3}

We may assume without loss of generality that $\ov{\Omega}$ lies in the half space $\{ x_1 >0 \}$.  Fix a unit vector $e$.  We wish to show that $u_{ee}\le 0$. 

  Define, for  a constant $A>0$,
$$Q_A = u_{ee} - A u- \sigma A x_1,$$
where $$\sigma= \frac{f(0)}{\sup_{x\in \ov{\Omega}} |x_1 f'(u(x))|+1}>0.$$

By our  assumptions $Q_A<0$ on $\partial \Omega$ for any $A>0$.  Choose $A$ sufficiently large, depending on $u$, so that $Q_A<0$ on $\Omega$.   We can do this because $x_1\ge c>0$ on $\ov{\Omega}$, for a uniform $c>0$.

Let's assume for a contradiction that if we let $A$ decrease towards zero, there is some $A=A_0>0$ at which $Q_A$ first attains 0 at some point in the interior.  Indeed if not then $Q_A<0$ for all $A>0$, and letting $A\rightarrow 0$ proves that $u_{ee} \le 0$,  as required. 

  Fix this $A=A_0$.  By assumption there is some interior point $p \in \Omega$ at which $0=Q_A(p)\ge Q_A(x)$ for all $x \in \Omega$.  We compute at this $p$, recalling the definition of $\sigma$, 
\[
\begin{split}
0 \ge {} & \Delta Q_A \\
 = {} & \Delta u_{ee} - A \Delta u \\
= {} & - f'(u) u_{ee} - f''(u) u_e^2 + A f(u) \\
\ge {} & - A f'(u) u - \sigma Af'(u) x_1 + Af(u) \\
> {} & A(- f'(u) u + f(u)) - A f(0) \ge0,
\end{split}
\]
a contradiction.  We used the concavity of $f$ twice.  First to observe that $-f''(u) u_e^2 \ge 0$ and second to see that 
$$\frac{f(u) - f(0)}{u} \ge f'(u).$$
which is the same as 
$$-uf'(u) + f(u) \ge f(0).$$
This completes the proof.

\section{A comparison to the contact parabola condition} \label{comparison}

Our proof of Theorem \ref{mainthm} is motivated by the proof of the main theorem in \cite{K}.  Kosmodem'yanskii showed that solutions to the torsion problem (\ref{torsion}) are  concave provided the following \emph{contact parabola condition} (CPC) holds:  for each point $p \in \partial \Omega$, the domain $\Omega$ is contained in the  \emph{third order contact parabola} at $p$.

 Working in preferred coordinates $x,y$ at $p$ (as in Definition \ref{preferred}), the third order contact parabola at $p$ is the curve with equation $y=(Ax+By)^2$ for constants $A, B$ with $A>0$ satisfying $y^{(k)}(0)=\rho^{(k)}(0)$ for $0\leq k \leq 3$.

   A key difference with condition (1) in Theorem \ref{mainthm}, abbreviated here by (CEC) (contact ellipse condition), is that (CEC) requires only first order contact between the ellipses with the boundary.   Indeed, in preferred coordinates around any $p\in \partial \Omega$ the ellipse in (CEC) will have the general equation \eqref{generalellipse} where the implicit function $y(x)$ must satisfy $y^{(k)}(0)=\rho^{(k)}(0)$ for $0\leq k \leq 1$ while the second and third derivatives are required only to be comparable, but not necessarily equal, through \eqref{assume1}.  
   
    To highlight the difference between these conditions, we will show that given any compact convex domain satisfying (CPC), an arbitrarily small $C^{3+\gamma}$ perturbation of $\Omega$, for $\gamma \in (0,1)$ may result in a domain which violates (CPC), while this is not the case for (CEC) (see Remark \ref{remark} (i)).

    Let $\Omega$ be a smooth compact convex domain satisfying (CPC).  We will construct a sequence of compact domains $\Omega_c$ such that for all $c>0$ sufficiently small, 
    \begin{enumerate}
    	\item[(i)] $\Omega_c$ is compact and convex and contains $\Omega$
    	\item[(ii)] $\Omega_c$ does not satisfy (CPC)
    	\item[(iii)] $\partial \Omega_c\to \partial \Omega$ uniformly in $C^{3+\gamma}$ as $c\to 0$, for any $\gamma \in (0,1)$.
    \end{enumerate}
    
     Let $p$ be any point on $\partial \Omega$, and  $y=(Ax+By)^2$ the 3rd order contact parabola at $p$  where $x, y$ are the preferred coordinates at $p$.   Writing this parabola near $x=0$ as $y=P(x)$ we compute, 
     $$P(x) = A^2 x^2 + 2A^3 Bx^3 +O(x^4).$$ 
If we write the boundary $\partial \Omega$ near the origin as $y=\rho(x)$, the assumption (CPC) implies that for $|x|\le c_0$, for a small constant $c_0>0$,
$$\rho(x) = A^2 x^2 + 2A^3B x^3 +  E(x),$$
where $E$ satisfies
$$|D^k E(x)| \le C x^{4-k},\quad k=0,1, \ldots, 4,$$
for a uniform constant $C$.
     
Define new smooth locally defined functions $\rho_c(x)$, for $0 < c \le c_0$, as follows.   Let $\zeta: \mathbb{R} \rightarrow [0,1]$ be a  smooth function equal to $1$ on $[-1/2,1/2]$ and equal to zero outside $[-1,1]$.  Then define for $|x| \le c_0$,
$$\rho_c(x) = \rho(x) - a \zeta(x/c)  x^4,$$
for a constant $a>0$ chosen sufficiently large so that 
\begin{equation} \label{violate}
\rho_c(x)=\rho(x) - ax^4 < P(x), \quad \textrm{for } |x| \le c/2.
\end{equation}

Note that $\rho_c \le \rho$ on $|x| \le c_0$  and $\rho_c = \rho$ outside $[-c,c]$.  Hence $\rho_c$ defines a new compact domain $\Omega_c$ which contains $\Omega$.

On $[-c,c]$ we have, for a uniform $C$,
\begin{equation} \label{rhoconvex}
\begin{split}
\rho_c''(x) = {} & \rho''(x) - a \left( \frac{\zeta''(x/c)}{c^2} x^4 - 12\zeta(x/c)x^2 - \frac{8 \zeta'(x/c)}{c} x^3 \right) \\
 \ge {} &  2A^2 - C c^2\ge A^2>0,
 \end{split}
 \end{equation}
shrinking $c_0$ if necessary.  Hence $\Omega_c$ is convex, and so (i) is satisfied.

The assertion (ii) follows from (\ref{violate}), which implies that $\Omega_c$ is not contained in the parabola $y=(Ax+By)^2$.

Finally, for (iii), we can estimate 
$$| \rho_c^{(k)}| \le C, \quad \textrm{for } k=0,1, \ldots, 4,$$
by a similar argument to that of (\ref{rhoconvex}) above.  Since $\rho_c \rightarrow \rho$ in $C^0([-c_0, c_0])$ as $c \rightarrow 0$, we have $\rho_c \rightarrow \rho$ in $C^{3 + \gamma}([-c_0, c_0])$ and this establishes (iii).

\begin{remark} \label{remark2}
We give here an explicit example of a domain which satisfies the assumption (CEC) of Theorem \ref{mainthm} but not (CPC).  Let  $\Omega$ be the unit disc, which satisfies (CPC) by \cite[Theorem 2]{K}.  Hence by the above argument and part (i) of Remark \ref{remark}, there is a small perturbation $\Omega_c$ of $\Omega$ which satisfies (CEC) but not (CPC).
\end{remark}


\begin{thebibliography}{0}
\bibitem{K} Kosmodem'yanskii, A.A., \emph{Sufficient conditions for the concavity of the solution of the Dirichlet problem for the equation $\Delta u=-1$}, Mat. Zametki 42 (1987), no. 4, 537--542.
\end{thebibliography}
\end{document}